# NONHOMOGENEOUS PARKING FUNCTIONS AND NONCROSSING PARTITIONS

DREW ARMSTRONG AND SEN-PENG EU

ABSTRACT. For each skew shape we define a nonhomogeneous symmetric function, generalizing a construction of Pak and Postnikov [9]. In two special cases, we show that the coefficients of this function when expanded in the complete homogeneous basis are given in terms of the (reduced) type of $k$-divisible noncrossing partitions. Our work extends Haiman's notion of a parking function symmetric function [5, 10].

## 1. INTRODUCTION

Let $\mu = (\mu_1, \mu_2, \ldots)$ and $\nu = (\nu_1, \nu_2, \ldots)$ be integer partitions; that is, weakly decreasing sequences of integers with only finitely many nonzero entries. We identify a partition $\mu$ with its Young diagram, an array of boxes with $\mu_i$ boxes in the $i$th row from the top, and each row aligned to the left.

Suppose that $\mu \subseteq \nu$ — this means that $\mu_i \leq \nu_i$ for all $i$, or that the Young diagram for $\mu$ is contained in that for $\nu$. The skew shape $\mu/\nu$ is the setwise difference of Young diagrams. A partial horizontal strip in the shape $\mu/\nu$ is a set of boxes, at most one in each column, such that the height of the boxes is weakly increasing to the right. We say a partial horizontal strip is right-aligned, and call it an r-strip, if by adding a box to the right of any given box, we no longer have a horizontal strip. If we think of the Young diagram as a grid of line segments, note that an r-strip is equivalent to a lattice path in the diagram from the bottom left corner to the top right (that is, the blocks of the r-strip indicate the horizontal steps of the path, see Figure 1). An ordinary horizontal strip is a partial horizontal strip with exactly one box in each column.

A block in an r-strip is a maximal sequence of adjacent boxes. The collection of sizes of blocks in an r-strip $\sigma$ form a partition, called the type of $\sigma$. By an abuse of notation, we will also denote this type by $\sigma$ when there is no confusion.

Let $\mathbf{x} = \{x_1, x_2, \ldots, \}$ be an infinite set of commuting variables. Given a skew shape $\mu/\nu$ we define a nonhomogeneous symmetric function

$$f_{\mu/\nu} = f_{\mu/\nu}(\mathbf{x}) := \sum_\sigma h_\sigma, \tag{1}$$

where the sum is over all r-strips $\sigma \subseteq \mu/\nu$ and $h_\sigma$ is the complete homogeneous symmetric function. For further definitions we refer to [11, Chapter 7].

Date: October 29, 2018.
D. Armstrong is partially supported by NSF grant DMS-0603567.
S.-P. Eu is partially supported by National Science Council, Taiwan under grants NSC 95-2115-M-390-006-MY3 and TJ & MY Foundation.

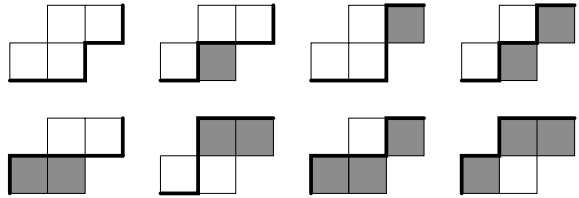

Figure 1. r-strips in the skew shape $(3,2)/(1)$.

For example, the symmetric function corresponding to the skew shape $(3,2)/(1)$ is given by
$$f_{(3,2)/(1)} = 2h_{(2,1)} + 2h_{(2)} + h_{(1,1)} + 2h_{(1)} + h_\emptyset.$$
This is clear from Figure 1, which displays the eight possible r-strips in this shape. We have also indicated the corresponding lattice paths by bold lines.

In this paper we will focus on two special shapes, namely $\mu/\nu = (n^{kn})/((n-1)^k, \ldots, 2^k, 1^k)$ — which we call the "stretched staircase" — and $\mu/\nu = (n^{kn})$. We call $f_{(n^{kn})/((n-1)^k,\ldots,2^k,1^k)}$ the Fuss type $A_{n-1}$ nonhomogeneous symmetric function and $f_{(n^{kn})}$ the Fuss type $B_n$ nonhomogeneous symmetric function. Our main results are to give formulas and a uniform explanation for the coefficients in these functions.

Our explanation will involve the noncrossing partitions. Given a partition of the set $[n] := \{1, 2, \ldots, n\}$, we say that the quadruple $(a, b, c, d)$, $1 \le a < b < c < d \le n$, is a "crossing" if $a, c$ are in one block and $b, d$ are in a different block of the partition. Let $\mathcal{NC}_n^A$ denote the set of noncrossing partitions of $[n]$ (that is, those without crossings). (Warning: There is an algebraic theory in which the set $\mathcal{NC}_n^A$ corresponds to the Cartan-Killing type $A_{n-1}$ — see [1]. We hope no confusion of indices will result.) A noncrossing partition can be represented pictorially by labelling the vertices of a convex $n$-gon clockwise by $1, 2, \ldots, n$ and associating to each block of a partition its convex hull. A partition is noncrossing if and only if its blocks are pairwise disjoint. The type of a noncrossing partition $\beta \in \mathcal{NC}_n^A$ is the integer partition given by its block sizes, and the reduced type is the type of the partition obtained from $\beta$ by deleting the block containing the symbol 1.

There is a 'Fuss generalization' for noncrossing partitions. Given a positive integer $k$, let $\mathcal{NC}_n^{A,(k)}$ denote the collection of noncrossing partitions of the set $[kn]$ in which each block has cardinality divisible by $k$. These are the $k$-divisible noncrossing partitions. The type of a $k$-divisible noncrossing partition $\beta \in \mathcal{NC}_n^{A,(k)}$ is the integer partition obtained by first taking the type of $\beta$ as an element of $\mathcal{NC}_{kn}^A$ and then dividing each entry by $k$. The reduced type again ignores the block containing symbol 1. For example, the 2-divisible noncrossing partition $1256/34/78$ has type $(2, 1, 1)$ and reduced type $(1, 1)$.

Given a partition $\lambda$, let $|\lambda|$ denote the sum of its parts, let $\ell(\lambda)$ denote its number of nonzero parts, and set $m_\lambda := m_1(\lambda)! m_2(\lambda)! m_3(\lambda)! \cdots$, where $m_i(\lambda)$ is the number of parts of $\lambda$ equal to $i$. Our main result on Fuss type $A_{n-1}$ symmetric functions is the following.

**Theorem 1.1.** *We have the expansion*

$$(2) \qquad f_{(n^{kn})/((n-1)^k,\ldots,2^k,1^k)} = \sum_{\lambda \vdash \le n} \frac{(k(n+1))! \cdot (n+1-|\lambda|)}{(n+1) \cdot m_\lambda \cdot (k(n+1) - \ell(\lambda))} h_\lambda.$$

*The sum is over all partitions $\lambda \vdash \leq n$, which means $\lambda \vdash n'$, for some $0 \leq n' \leq n$. Moreover, the coefficient of $h_\lambda$ is the number of noncrossing partitions in $\mathcal{NC}_{n+1}^{A,(k)}$ with reduced type $\lambda$.*

Note that the *reduced* type of a noncrossing partition of the set $[n+1]$ is an integer partition $\lambda$ of $n'$ for some $n' \leq n$. Thus the sum of the coefficients in (2) is the Catalan number $\frac{1}{n+2}\binom{2(n+1)}{n+1}$.

Type $B$ noncrossing partitions were introduced by Reiner (see [1]), and they can be represented in the following way. Consider a regular $2n$-gon with vertices labeled by $1, 2, \ldots, n, -1, -2, \ldots, -n$ clockwise. A **type $B$ noncrossing partition** is a noncrossing partition on these vertices that is invariant under the antipodal map. Let $\mathcal{NC}_n^B$ denote the collection of noncrossing partitions of the set $[n]^\pm := \{-n, \ldots, -2, -1, 1, 2 \ldots, n\}$. The antipodal map decomposes the blocks into orbits of size two, together with possibly one orbit of size one whose element we call the **antipodal block**. The **type** of a type $B$ noncrossing partition is calculated as in type $A$ by discarding the antipodal block and counting only one block from each orbit of size two. Note that this type is in some sense already "reduced".

The type $B$ noncrossing partitions also have a Fuss generalization (see [1]). Given a positive integer $k$, a **k-divisible type $B$ noncrossing partition** is a type $B$ noncrossing partition in which each block has cardinality divisible by $k$. Let $\mathcal{NC}_n^{B,(k)}$ denote the set of $k$-divisible type $B$ noncrossing partitions of the set $[kn]^\pm$. The **type** of a $k$-divisible type $B$ noncrossing partition is calculated by first considering its type as an element of $\mathcal{NC}_{kn}^B$ and then dividing each entry by $k$.

Our main result on Fuss type $B_n$ symmetric functions is as follows.

**Theorem 1.2.** *We have the expansion*

$$(3) \qquad f_{(n^{kn})} = \sum_{\lambda \vdash \leq n} \frac{(kn)!}{m_\lambda (kn - \ell(\lambda))!} h_\lambda.$$

*The sum is over all partitions $\lambda \vdash \leq n$, which means $\lambda \vdash n'$, for some $0 \leq n' \leq n$. Moreover, the coefficient of $h_\lambda$ is the number of noncrossing partitions in $\mathcal{NC}_n^{B,(k)}$ with type $\lambda$.*

The rest of this paper is organized as follows. The background is put in Section 2 and we will see that this work is a generalization of Haiman's and Stanley's results. In Sections 3 and 4 we prove the Fuss type $A$ and Fuss type $B$ cases respectively.

## 2. Some Background

### 2.1. Parking functions.

A sequence of positive integers $\alpha = (a_1, \ldots, a_n)$ is called a **parking function** of length $n$ if its nondecreasing rearrangement $(b_1, \ldots, b_n)$ satisfies $b_i \leq i$ for all $i$. We denote by $\mathcal{P}_n$ the set of parking functions of length $n$. If $m_i$ denotes the number of entries in $\alpha$ equal to $i$, then the **type** of the parking function $\alpha$ is the integer partition given by the weakly decreasing rearrangement of the $m_i$'s. Parking functions were first introduced by Konheim and Weiss [7].

A parking function is called **primitive** if the $a_i$'s are weakly increasing. Hence a parking function is a permutation of a primitive parking function. It is well known that the number of parking functions of length $n$ is $(n+1)^{n-1}$ and the number of primitive parking functions of length $n$ is the Catalan number $\frac{1}{n+1}\binom{2n}{n}$, see [10].

A primitive parking function of length $n$ can be seen as a choice of a set of boxes, one from each column, of weakly increasing heights in the staircase shape $(n^n)/((n-1), \ldots, 2, 1)$. Hence it is equivalent to a horizontal strip.

From this viewpoint we may define the set of parking functions with respect to a skew shape $\mu/\nu$. A primitive parking function with respect to $\mu/\nu$ is the sequence of (weakly increasing) heights of the boxes in a horizontal strip, and a parking function is a permutation of a primitive one. The $(n^n)$-parking functions are called "type $B$ parking functions" in [4]. The enumeration of $\mu/\nu$ parking functions was done by Kung et al. [8].

2.2. **The parking function symmetric function.** Consider the set of parking functions of length $n$ together with the action of the symmetric group $\mathfrak{S}_n$ by permuting entries. We denote by $pf_n(\mathbf{x})$ the Frobenius chracteristic of the action, and call it the parking function symmetric function. It was first considered by Haiman [5]. Explicit expansions of $pf_n(\mathbf{x})$ into the classical bases were given by Stanley. Among them was the following, together with a combinatorial interpretation of the coefficients.

**Theorem 2.1.** [10, Proposition 2.4] *We have the expansion*

$$(4) \qquad pf_n(\mathbf{x}) = \sum_{\lambda \vdash n} \frac{n!}{m_\lambda(n+1-\ell(\lambda))} h_\lambda.$$

*The sum is over all partitions $\lambda$ of $n$. Moreover, the coefficient of $h_\lambda$ is the number of noncrossing partitions $\in \mathcal{NC}_n^A$ of type $\lambda$.*

If $\chi$ is the permutation character of $\mathfrak{S}_n$ acting on the parking functions of length $n$, then we have $\chi = \sum \chi_\mathcal{O}$, where the sum is over $\mathfrak{S}_n$-orbits, and $\chi_\mathcal{O}$ is the character of $\mathfrak{S}_n$ acting on orbit $\mathcal{O}$. Note that the parking functions in an orbit all have the same type. If the parking functions in orbit $\mathcal{O}$ have type $\lambda$, Stanley observed ([10, Proposition 2.4]) that the Frobenius characteristic of $\chi_\mathcal{O}$ is the complete homogeneous symmetric function $h_\lambda$.

Since there is a unique primitive parking function in each $\mathfrak{S}_n$-orbit, we conclude that the coefficient of $h_\lambda$ in (4) counts the number of primitive parking functions of type $\lambda$. Finally, there is a well known bijection from primitive parking functions of type $\lambda$ to noncrossing partitions of type $\lambda$ (see Section 3 below).

If $\lambda$ is a partition of $n$, note that the number of noncrossing partitions $\in \mathcal{NC}_n^A$ of type $\lambda$ is equal to the number of noncrossing partitions $\in \mathcal{NC}_{n+1}^A$ of *reduced* type $\lambda$ (to each of the former, we add the singleton block $\{1\}$). Thus, comparing with our Theorem 1.1, we conclude that Haiman's function $pf_n(\mathbf{x})$ is the top homogeneous part of our $f_{\mu/\nu}$ when $\mu/\nu$ is the staircase shape $(n^n)/((n-1),\ldots,2,1)$. Stanley also discussed a generalization which is the top homogeneous part of our $f_{\mu/\nu}$ when $\mu/\nu$ is the stretched staircase [10, Section 5].

3. Fuss type $A_{n-1}$

In this section we prove Theorem 1.1. Let $\mathcal{H}_n^{(k)}$ denote the set of all r-strips in the stretched staircase $(n^{kn})/((n-1)^k,\ldots,2^k,1^k)$.

The proof involves Fuss-Catalan paths. Given a positive integer $k$, we say that a $k$-Fuss-Catalan path of length $n$ is a lattice path from $(0,0)$ to $(n,kn)$, using east steps $E = (0,1)$ and north steps $N = (1,0)$, and staying in the region $0 \le y \le kx$. We denote by $\mathcal{D}_n^{(k)}$ the set of $k$-Fuss-Catalan paths of length $n$. A maximal sequence of consecutive east steps is called an ascent. The collection of lengths of all ascents defines the type of the path. Note that the first step of the path must always be $E$. To compute the reduced type of the path, we ignore the ascent containing this first $E$. (When we present a bijection to the noncrossing partitions, this first east step will be labelled by the symbol 1). When $k = 1$, a Fuss-Catalan paths is called a Dyck path.

The proof will consist of three steps. First we will biject the set of r-strips $\mathcal{H}_n^{(k)}$ with the set $\mathcal{D}_{n+1}^{(k)}$ of $k$-Fuss-Catalan paths of length $n+1$. Second, we will give a bijection of $\mathcal{D}_{n+1}^{(k)}$ with the $k$-divisible noncrossing partitions $\mathcal{NC}_{n+1}^{A,(k)}$, preserving both type and reduced type. Third, we count the noncrossing partitions $\mathcal{NC}_{n+1}^{A,(k)}$ with respect to reduced type. Note the transition of indices in the first bijection.

**Lemma 3.1.** *There is a bijection $\phi_n^{(k)} : \mathcal{H}_n^{(k)} \to \mathcal{D}_{n+1}^{(k)}$, which sends the type of the r-strip $\nu \in \mathcal{H}_n^{(k)}$ to the reduced type of of the path $\alpha := \phi_n^{(k)}(\nu) \in \mathcal{D}_{n+1}^{(k)}$.*

*Proof.* Fix the bottom left and upper right corners of the skew shape $(n^{kn})/((n-1)^k, \ldots, 2^k, 1^k)$ at $(1, 0)$ and $(n+1, kn)$ in $\mathbb{Z}^2$, respectively, and add the two line segments $(0,0) \to (1,0)$ and $(n+1, kn) \to (n+1, k(n+1))$. Any lattice path from $(0,0)$ to $(n+1, k(n+1))$ on this configuration is a $k$-Fuss-Catalan path $\in \mathcal{D}_{n+1}^{(k)}$. Denote the blocks of the given r-strip $\nu \in \mathcal{H}_n^{(k)}$ by $B_1, B_2, \ldots, B_s$ from left to right. Then there is a unique lattice path $\alpha := \phi_n^{(k)}(\nu)$ from $(0,0)$ to $(n+1, k(n+1))$ traveling along the lines $y = 0, L_W(B_1), L_N(B_1), \ldots, L_W(B_k)$, $L_N(B_s)$, $x = n+1$ in order, where $L_W(B_i)$ and $L_N(B_i)$ are the lines passing the left and north boundary of $B_i$, $1 \leq i \leq k$, respectively. The map is clearly invertible. $\square$

Compare with the example in Figure 1.

**Lemma 3.2** (Labeling Lemma). *There is a bijection $\psi_n^{(k)} : \mathcal{D}_n^{(k)} \to \mathcal{NC}_n^{A,(k)}$ which simultaneously preserves type and reduced type.*

*Proof.* Given $\alpha \in \mathcal{D}_n^{(s)}$. Note that the plane is divided into disjoint *regions* by the lines $y = kx + i$, $i \in \mathbb{Z}$ and each east step is cut by these lines into $k$ segments, each of which is in a region and is of length $\frac{1}{k}$. We put a vertex centered in every segment and add edges connecting these vertices to obtain a rooted tree, called the auxiliary tree, by the following rules:

(1) Each vertex is connected to adjacent vertices in the same ascent.
(2) The leftmost vertex of each ascent connects with the first lower vertex (if exists) in the same region.

The resulting tree, rooted at the vertex centered on the first segment $(0,0) \to (\frac{1}{k}, 0)$ is a binary tree. We label the root by the symbol 1 and other vertices by $\{2, \ldots, k(n+1)\}$ in preorder, or depth-first search (i.e., always choose the upward branch if exists and proceed recursively.) Now we define a partition $\beta = \psi_n^{(k)}(\alpha)$ whose blocks are the labels on the ascents of $\alpha$. The partition $\beta$ is $k$-divisible because the cardinality of each block is a multiple of $k$ and is noncrossing because of preorder. Hence $\beta \in \mathcal{NC}_n^{A,(k)}$

The inverse mapping $(\psi_n^{(k)})^{-1}$ can be constructed as follows. Given $\beta \in \mathcal{NC}_n^{A,(k)}$, one can list the blocks $B_1, B_2, \ldots, B_s$ in a line such that $\min B_i < \min B_{i+1}$ for all $i$, and such that the numbers within each block are listed in increasing order. We say that such a listing is canonical.

Suppose $\beta$ is in canonical order. Let $a_i$ and $b_i$ denote the smallest and the largest entries, respectively, of the block $B_i$. For each integer $1 \leq t \leq kn$ we create a vertex $P_t := (x_t, y_t)$ in $\mathbb{Z}^2$, and we connect these vertices as follows to obtain the auxiliary tree:

(1) Let $P_1 := (0,0)$ be the root. Other numbers in $B_1$ correspond to the vertices to the right of $P_1$ one by one on the line $y = 0$ in order, each connecting to the vertex to its left by an edge of unit length.

(2) For $2 \leq i \leq k$, let $P_{a_i} := (x_{a_i-1}, s+1-i)$, other numbers in $B_i$ correpond to the vertices to the right of $P_{a_i}$ one by one on the line $y = s+1-i$ in order, each connecting to the vertex to its left by an edge of unit length.

(3) For $2 \leq i \leq k$, connect vertices $a_i$ and $a_i - 1$ by an edge.

The binary tree obtained is the auxilary tree, rooted at 1 and labeled in preorder. Setting $x_{a_{s+1}} = 0$, the desired path $(\psi_n^{(k)})^{-1}(\beta)$ is constructed by drawing $\frac{1}{k}|B_i|$ east steps followed by $x_{b_i} - x_{a_{i+1}} + 1$ north steps, $1 \leq i \leq s$, in order.

Since the first segment of the first step of $\alpha$ is always labeled 1 and corresponds to the block containing 1 in $\beta$ and vice versa, $\psi_n^{(k)}$ preserves both type and reduced type simultaneously. $\square$

For example, Figure 2 illustrates how the path

$$ENEENNNNENNNEENNNN \in \mathcal{D}_6^{(2)}$$

on the left bijects with the partition $1, 6/2, 3, 4, 5/7, 10, 11, 12/8, 9 \in \mathcal{NC}_6^{A,(2)}$ on the right, with the auxilary tree displayed in the middle. Note that the type is $(2, 2, 1, 1)$ and the reduced type is $(2, 2, 1)$.

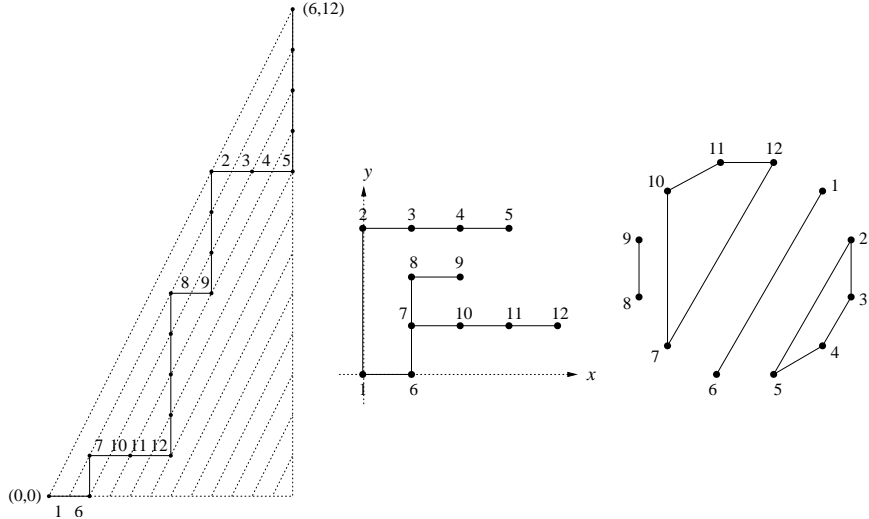

FIGURE 2. An example of the type A bijection.

In the following lemma we count the $k$-divisible noncrossing partitions by reduced type. Recall that the reduced type is obtained by discarding the block containing the symbol 1 (or any one arbitrary fixed symbol).

**Lemma 3.3.** *Let $p^{(k)}(\lambda)$ denote the number of $k$-divisible noncrossing partitions $\in \mathcal{NC}_n^{A,(k)}$ with reduced type $\lambda$. Then we have*

$$p^{(k)}(\lambda) = \frac{(kn)! \cdot (n - |\lambda|)}{n \cdot m_\lambda \cdot (kn - \ell(\lambda))},$$

*where $|\lambda|$ is the sum of the parts of $\lambda$, $\ell(\lambda)$ is the number of nonzero parts of $\lambda$, and $m_\lambda = m_1(\lambda)! m_2(\lambda)! \cdots$, where $m_i(\lambda)$ is the number of times that $i$ occurs in $\lambda$.*

*Proof.* If the $k$-divisible noncrossing partition $\pi \in \mathcal{NC}_n^{A,(k)}$ has reduced type $\lambda$, let $\zeta \vdash n$ be the type of $\pi$. This is the partition obtained from $\lambda$ by adding one part of size $k(n - |\lambda|)$. Thus $\ell(\zeta) = \ell(\lambda) + 1$ and $m_\zeta = m_{n-|\lambda|}^\zeta m_\lambda$, where $m_{n-|\lambda|}^\zeta$ is the number of times that part $n - |\lambda|$ appears in $\zeta$.

By [1, Theorem 4.4.4], which follows from [6, Theorem 4], the number of $k$-divisible noncrossing paritions with type $\zeta$ is equal to
$$q^{(k)}(\zeta) = \frac{(kn)!}{m_\zeta(kn + 1 - \ell(\zeta))!}.$$

Thus it suffices to show that
$$kn \cdot p^{(k)}(\lambda) = q^{(k)}(\zeta) \cdot m_{n-|\lambda|}^\zeta \cdot k(n - |\lambda|).$$

We prove this by a double counting argument. On one hand, $kn \cdot p^{(k)}(\lambda)$ is the number of *pointed* noncrossing partitions of $[kn]$ with reduced type $\lambda$ (that is, the number of partitions in which we have some distinguished symbol $i \in [kn]$ and the blocks not containing $i$ have type $\lambda$). On the other hand we may create such a partition by starting with an arbitrary partition of type $\zeta$, of which there are $q^{(k)}(\zeta)$. Then we select a block of size $k(n - |\lambda|)$ in $m_{n-|\lambda|}^\zeta$ ways. From within this block we may choose a distinguished symbol in $k(n - |\lambda|)$ ways. Thus the number of pointed noncrossing partitions with reduced type $\lambda$ is also equal to $q^{(k)}(\zeta) \cdot m_{n-|\lambda|}^\zeta \cdot k(n - |\lambda|)$. □

*Proof. of Theorem 1.1*: By Lemma 3.1 and Lemma 3.2, the map $\psi_{n+1}^{(s)} \circ \phi_n^{(s)}$ is a bijection $\mathcal{H}_n^{(s)} \to \mathcal{NC}_{n+1}^{A,(k)}$, sending type into reduced type. Now apply Lemma 3.3. □

## 4. Fuss type $B_n$

In this section we prove Theorem 1.2. The proof involves Fuss binomial paths. Given positive integer $k$, a $k$-**Fuss binomial path of length** $(k + 1)n$ is a lattice path from $(0, 0)$ to $(n, kn)$, using east steps $E = (1, 0)$ and north steps $N = (0, 1)$. Denote by $\mathcal{B}_n^{(k)}$ the set of $k$-Fuss bimonial paths of length $(k + 1)n$.

A maximal sequence of consecutive east steps is called an **ascent**. The collection of the lengths of all ascents except for the possible one on $y = 0$ is again a partition, which we define to be the **type** of this path. Note that it is in fact the "reduced type" if one wants to be consistent with the type $A_{n-1}$ case. The reason for calling it 'type' rather than 'reduced type' is to be in consistent with the definition of 'type' in $\mathcal{NC}_n^B$.

The structure of the proof is similar to that of the Fuss type $A_{n-1}$ case. The proof of the following lemma is omitted.

**Lemma 4.1.** *There is a type preserving bijection $\phi_n^{(k)} : \mathcal{H}_n^{(k)} \to \mathcal{B}_n^{(k)}$.*

The following proposition is central to the proof.

**Lemma 4.2.** *There is a bijection $\psi_n^{(k)} : \mathcal{B}_n^{(k)} \to \mathcal{NC}_n^{B,(k)}$ which preserves type.*

*Proof.* Consider a path $\alpha \in \mathcal{B}_n^{(k)}$. The plane is divided into disjoint regions by the lines $y = skx + i$, $i \in \mathbb{Z}$ and each east step of $\alpha$ is cut by these lines into $s$ segments, each in a region. We call the intersection of the rectangle $[0, n] \times [0, kn]$ with $0 \leq y \leq kx$ the **positive triangle**, and its intersection with $kx \leq y$ the **negative triangle**. Note that $\alpha$ is decomposed into a concatenation of positive and negative paths, where now a segment of length $\frac{1}{k}$ serves as an east step. By a positive (negative) path we mean that all of its east steps are in the postive (negative) triangle. Clearly the number of positive and negative paths differs by at most 1, and we may write $\alpha = P_1 N_1 P_2 N_2 \ldots P_k N_s$ for some $s$, where $P_i$ (or $N_i$), $1 \leq i \leq k$, are positive (negative). Note that both, one, or neither of $P_1$ and $N_s$ might be empty.

Each $P_i$ (or $N_i$) is a Dyck path (rotated by 180 degrees in the case of $N_i$), using an *east segment* of length $\frac{1}{k}$ as the unit of an east step and 1 as the unit of a north step. Denote by $p_i$ (or $n_i$) the number of east segments in $P_i$ (or $N_i$) and set $n_0 := n_1 + \cdots + n_s$. For each $P_i$ and $N_i$, seen as a Dyck path, we apply the Labeling Lemma 3.2 to label east segments: we label the east segments of $N_s$ by $[1, n_s]$; of of $N_{s-i}$ by $[n_s + n_{s-1} + \cdots + n_{s-i+1} + 1, n_s + n_{s-1} + \cdots + n_{s-i+1} + n_{s-i}]$ for $1 \leq i \leq s-1$; of $P_1$ by $[n_0 + 1, n_0 + p_1]$; and of $P_i$ by $[n_0 + p_1 + \cdots + p_{i-1} + 1, n_0 + p_1 + \cdots + p_{i-1} + p_i]$ for $2 \leq i \leq k$.

Now we define a partition $\beta = \psi_n(\alpha) \in \mathcal{NC}_n^{B,(k)}$ by the following rules:

(1) The ascent on $y = 0$ corresponds to the antipodal block (the unique block containing both negative and positive integers), if it exists.
(2) Each ascent not on $y = 0$ corresponds to a pair of symmetric blocks, one of which consists of those numbers labelling this ascent.

It is clear that type is preserved under $\psi_n^{(k)}$. And $\beta$ is noncrossing because the subpartition using positive (negative) numbers is noncrossing, and because of the fact that the positive and negative paths are alternating.

Now we consider the inverse mapping $(\psi_n^{(k)})^{-1}$. Below we will use the ordering on indices $-1 < -2 < \cdots < -kn < 1 < 2 < \cdots < kn$. Given $\beta \in \mathcal{NC}_n^{B,(k)}$, we list the blocks $B_1, \ldots, B_s$ of $\beta$ such that $\min B_i < \min B_{i+1}$, and such that the numbers in each block are listed in clockwise order with respect to their positions as vertices of the regular $2kn$-gon. We call such a listing of $\beta$ **canonical**. For example, $-1, 6, 9, 10/-2, -3/-4, -5, 4, 5/6, 9, -10, 1/-7, -8/2, 3/7, 8 \in \mathcal{NC}_5^{B,(2)}$ is a canonical listing.

Now suppose $\beta$ is canonical. We will select one element from each antipodal pair $\{i, -i\}$, $1 \leq i \leq kn$, serving as the labels on the east segments of $(\psi_n^{(k)})^{-1}(\beta)$. The idea is to collect numbers block by block, starting from the block containing $-1$, until we have half of them.

(1) If there is an antipodal block, let $a_0$ be the smallest positive number in the antipodal block.
(2) If there is no antipodal block, find the last block containing both positive and negative numbers such that the absolute value of any negative number in this block is smaller than any positive number in this block. Let $a_0$ be the smallest positive number in this block.
(3) If there is no block containing both positive and negative numbers, set $a_0 = n + 1$.

For $a < b$, let $[a, b]$ denote the set of integers $a \leq x \leq b$. We define $A := [-a_0 + 1, -1] \cup [a_0, n]$ to be the set of selected numbers. Note that by definition the positive (negative) numbers in $A$ also form a noncrossing partition, and can be listed in canonical order. Now we put $kn$ numbers in $A$ into a $2 \times kn$ array, leaving half of the positions empty, subject to the following rules:

(1) The negative numbers, starting from the smallest one, are put in the row 1 from right to left in the canonical order.
(2) The positive numbers, starting from the smallest one, are put in the row 2 from left to right in the canonical order.

(3) Both positions $(1,j)$ and $(2,j+1)$ are occupied if and only if the following condition holds: $(1,j)$ is the smallest negative number of some block containing simultaneously postive and negative numbers, and $(2,j+1)$ is the smallest positive number in this block.

In this way the $2 \times n$ array is a concatenation of sections of numbers, each is a maximal sequence of numbers of the same sign, and the signs are alternating among sections. Let us call the sections from left to right $P_1, N_1, P_2, N_2, \ldots, P_s, N_k$, where $P_i$ (or $N_i$) stands for the $i$-th positive (negative) section. Note that both, one or neither of $P_1, N_k$ might be empty.

For each section $P_i$ (or $N_i$), using the Labelling Lemma 3.2 we can construct a Dyck path with east segments labeled. Rotate the path constructed from $N_i$ by 180 degrees. By abuse of notation we also call the rotated path $N_i$. Our desired $\alpha = (\psi_n^{(k)})^{-1}(\beta) \in \mathcal{B}_n$ is defined to be the concatenation of the paths $\alpha := P_1 N_1 P_2 N_2 \ldots P_s N_s$. It is easy to see that this construction preserves type. $\square$

For example, let $k = 2$. Figure 3 illustrates how the path

$$ENNNNNNENNENNNEN \in \mathcal{B}_4^{(3)}$$

bijects with the 3-divisible noncrossing partition

$$-1, -2, 12/-3, -7, 11/-4, -5, -6/-8, -9, -10, 8, 9, 10/-11, 3, 7/-12, 1, 2/4, 5, 6 \in \mathcal{NC}_4^{B,(3)}$$

with the auxillary $2 \times 12$ array in the middle.

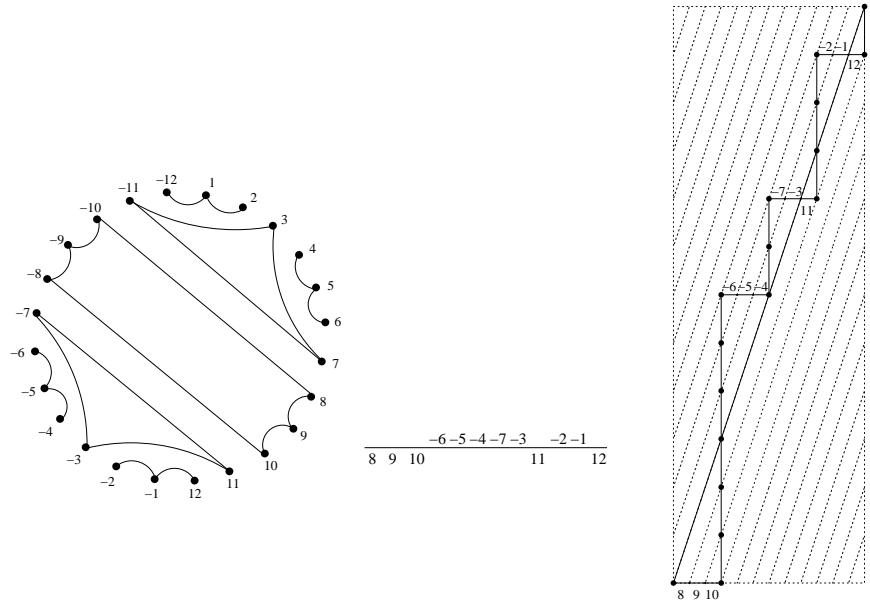

FIGURE 3. An example of the type B bijection.

Now it suffices to count number of $k$-divisible type $B$ noncrossing partitions with respect to type. Luckily this follows from a result of Athanasiadis (see [1, Theorem 4.5.11])

**Lemma 4.3.** [2] *The number of $k$-disible type $B$ noncrossing partitions $\in \mathcal{NC}_n^{B,(s)}$ with type $\lambda$ is*

$$\frac{(kn)!}{m_\lambda(kn-\ell(\lambda))!}. \tag{5}$$

*Proof. of Theorem 1.2*: By Lemma 4.1 and Lemma 4.2, the map $\psi_n^{(k)} \circ \phi_n^{(k)}$ is a bijection from $\mathcal{H}_n^{(k)}$ to $\mathcal{NC}_n^{B,(k)}$ preserving type. Now apply Lemma 4.3 and we are done. $\square$

## 5. Concluding remarks

As we can see, the 'Fuss type $A$' shape $(n^{kn})/((n-1)^k, \ldots, 2^k, 1^k)$ is the dilation of the 'type $A$' shape $(n^n)/(n-1, \ldots, 2, 1)$, and the 'Fuss type $B$' shape $(n^{kn})$ is the dilation of the 'type $B$' shape $(n^n)$ by a factor of $k$ along the $y$ direction. It is interesting to discuss the effect on $f_{\zeta/\mu}$ when applying various operations on the skew shape $\zeta/\mu$. Also it is interesting to investigate $f_{(\zeta/\mu)\star(\zeta'/\mu')}$ where $\star$ is some binary operation. Although a general theory can be developed, in this work we only focus on two special shapes owing to their direct connections to the noncrossing partitions.

Stanley [10, Theorem 2.3] also considered the relationship of Haiman's parking function symmetric function to the Ehrenborg quasisymmetric function of the lattice of noncrossing partitions. There he noticed a curious discrepancy in the indices. We hope that our nonhomogeneous version of the parking functions can help to explain this discrepancy. Perhaps there is a nonhomogeneous version of the Ehrenborg function.

## Acknowledgements.

This paper was prepared during the second named author's visit to the School of Mathematics, University of Minnesota. The authors thank Victor Reiner for helpful suggestions.

School of Mathematics, University of Minnesota, Minneapolis, Minnesota, USA
*E-mail address*: `armstron@math.umn.edu`

Department of Applied Mathematics, National University of Kaohsiung, Kaohsiung 811, Taiwan, ROC
*E-mail address*: `speu@nuk.edu.tw`